\documentclass[11pt,oneside,psamsfonts,final]{amsart}

\usepackage[body={5.45in,9.3in}]{geometry}

\usepackage[mathscr]{eucal}
\usepackage{amssymb,amsmath,amsthm,graphicx,times,mathptm}

\usepackage{amsfonts}

\usepackage{hyperref}
\usepackage[all]{hypcap} 

\usepackage[oldenum]{paralist} 

\usepackage{soul} 
\sodef\TopTiFo{}{0.03em}{.5em}{2em plus 0.1em minus 0.1em}

\usepackage{hyperref}

\DeclareMathAccent{\Interior}{\mathord}{operators}{'027}

\makeatletter

\def\@settitle{\begin{center}%
\uppercasenonmath\@title{
  \large\textrm{\@title}}
  \end{center}%
}
\copyrightinfo{}{}

\def\@setauthors{%
  \begingroup
  \trivlist
  \centering\footnotesize \@topsep40\p@\relax
  \advance\@topsep by -\baselineskip
  \item\relax
  \andify\authors
\textsc{\authors}%
  \endtrivlist
{\hfill \scriptsize 
        (\textit{Received $16$ February $1981$}; \hfill}
\par
{\hfill \scriptsize 
        \textit{annotated and \textup{arXiv}ed by the author %
               $13$ November $2004$}) \hfill}
  \endgroup
}
\def\@seccntformat#1{\bfseries
       \protect\footnotesize{$\mathbf{\S}$\csname the#1\endcsname.\quad}}
\newdimen\@bls                              
\@bls=\baselineskip                         
\def\section{\@startsection{section}{1}{\z@}{1.5\@bls\@plus .4\@bls 
             \@minus .1\@bls}{\@bls}{\centering\footnotesize\bfseries}}

\def\@captionheadfont{\normalfont\footnotesize}
\def\@captionfont{\normalfont\footnotesize}
\newtheoremstyle{TopologyThm}{\baselineskip}{\baselineskip}{\itshape}%
  {\parindent}{\fontshape{sc}\selectfont}{.}{ }{}
\newtheoremstyle{TopologyRmk}{\baselineskip}{\baselineskip}{}%
  {\parindent}{\fontshape{it}\selectfont}{.}{ }{}

\theoremstyle{TopologyThm}
\newtheorem*{theorem*}{T{\small{heorem}}}

\newtheorem{lemma}{L{\small{emma}}}

\theoremstyle{TopologyRmk}

\newtheorem*{remarks}{Remarks}
\newtheorem{example}{Example}[section]

\def\footnoterule{\kern-.4\p@
        \hrule\@width 5.5pc\kern11\p@\kern-\footnotesep}
\def\@setthanks{\def\thanks##1{\kern-\parindent$\dagger$##1\@addpunct.}\thankses}

\def\article@logo{%
  \set@logo{%
    \textup{Hyper{\TeX}ed arXival version prepared and 
    updated December 2003 and November 2004.}\newline
    \textup{Originally published in \textit{Topology}, 
     Vol.\ 22, No.\ 2, pp.\ 191--202, 1983.}\newline
    \textup{Marginal notes refer to the 
    \protect\hyperlink{addenda}{Addenda}.}
  }%
}

\DeclareRobustCommand{\appnote}[1]{%
\marginpar{\quad\label{appnote:#1}%
\hypertarget{appnote:#1}{\large\textcircled{\small\ref{#1}}}}}

\def\qedsymbol{{\ensuremath\blacksquare}}
\newenvironment{unproved}[1]{\def\QQQQ{\string#1}\begin{\QQQQ}}%
                                  {\hfill\qedsymbol\end{\QQQQ}}

\renewcommand\a{\alpha}
\newcommand\Bd{\partial} 
\newcommand{\deffont}{\textit}        
\newcommand{\bydef}[1]{\deffont{#1}}  
\newcommand{\Suchthat}{\negthinspace\colon} 
\newcommand{\C}{\mathbb{C}}
\def\s#1{{\sigma\mathstrut}_{#1}}
\newcommand{\e}{\varepsilon}
\newcommand{\sub}{\subset}
\newcommand{\R}{\mathbb{R}}
\renewcommand\P{\mathbb{P}}
\newcommand{\sgn}{\operatorname{sgn}}
\newcommand{\pr}{\operatorname{pr}}

\begin{document}

\bibliographystyle{amsplain}
\title{\TopTiFo{Algebraic functions and closed braids}}
\author[Lee Rudolph]%
{\large{L}\small{ee} \large{R}\small{udolph}%
\thinspace$\dagger\ddagger$}

\thanks{Partially supported by NSF grant 
{MCS 76-08230}.\\
$\ddagger$\thinspace Present address: Box 251, Adamsville,
RI 02801, U.S.A}

\def\ufootnote#1{\let\savedthfn\thefootnote\let\thefootnote\relax
\footnote{#1}\let\footnotemark\savedthfn\addtocounter{footnote}{-1}}

\def\dots{\mathinner{\mkern2mu\ldotp\mkern2mu\ldotp\mkern2mu\ldotp\mkern2mu}}
\def\dotsm{\mathinner{\mkern2mu\cdotp\mkern2mu\cdotp\mkern2mu\cdotp\mkern2mu}}

\maketitle
\makeatother

\setcounter{section}{0}

\section{INTRODUCTION}

\noindent{\textsc{\large{L}\small{et}}} 
$f(z,w)\equiv f_0(z)w^n + f_1(z)w^{n-1} + \dotsm + f_n(z)\in \C[z,w]$.  
Classically, the equation $f(z,w)=0$ was said to define 
$w$ as an ($n$-valued) \bydef{algebraic function} of $z$, 
provided that
$f_0(z)$ was not identically $0$ and that $f(z,w)$ 
was squarefree and without factors of the
form $z-c$.  Then, indeed, the \bydef{singular set}
$B=\{z\Suchthat \text{there are not $n$ distinct 
solutions $w$ to}$$f(z,w)=0\}$ is finite; 
and as $z$ varies in any simply-connected
domain avoiding $B$, the
$n$ distinct solutions $w_1,\dots,w_n$ of
$f(z,w)=0$ will be analytic functions of $z$.  
Now let $\gamma$
be a simple closed curve in $\C-B$.  In the open solid torus
$\gamma\times \C\sub\C^2$, the set $K_\gamma=V_f\cap\gamma\times\C$
(where $V_f=\{(z,w)\Suchthat f(z,w)=0\}$) is evidently a closed
$1$-manifold, as smooth as $\gamma$, on which the projection to
$\gamma$ is an $n$-sheeted (possibly disconnected) covering map.
A $1$-manifold in a solid torus, which projects as a covering onto
the circle factor, is called a \bydef{closed braid}.  When the torus
is embedded (in the standard way) in a $3$-sphere (as $\gamma\times \C$
will be, shortly), the closed braid becomes a knot or link in that 
sphere; if the circle factor is oriented, there is a natural way to
orient that knot or link.  Which such oriented links, we may ask, 
arise from algebraic functions (when $\gamma$ is oriented counterclockwise)?

The points $z_0\in B$ are of two kinds (some may be of both).
If, for some $w_0$ such that $f(z_0, w_0) = 0$, 
it also happens that $(\partial f/\partial w)(z_0,w_0)$, 
we call $z_0$ a \bydef{singular point} 
of the algebraic function. (Either $(z_0, w_0)$ is a 
singular point, in the usual sense, of the algebraic curve $V_f$,
or it is a regular point at which the tangent line is the vertical 
line $z=z_0$.)  At a singular point $z_0$, some solution $w$ 
to $f(z_0, w) = 0$ has multiplicity greater than $1$. 
On the other hand, $z_0$ may be a root of $f_0(z)$; 
then there are not $n$ solutions, even counting 
multiplicities, to $f(z_0, w) = 0$. 
A root of $f_0(z)$ is a \bydef{pole} of the algebraic
function.

The set $K_\gamma$, being compact, actually lies in 
some closed solid torus 
$\gamma \times D_r = \{(z, w)\Suchthat z\in \gamma, |w|\le r\}$. 
Let $B^4$ be the bicylinder $D \times D_r$ where $D$ is the bounded
region in $\C$ with $\Bd D = \gamma$; then $B^4$ is 
homeomorphic to a $4$-ball, and its boundary $3$-sphere is 
decomposed in the usual way into two solid tori, 
$\gamma \times D_r$ and $D \times \Bd D_r$. 
\textit{If no pole of $f(z,w)$ lies in $D$, 
then $K_\gamma$ is the entire intersection of $V_f$ with 
$\Bd B$}; that is, $V_f$ does not meet $D \times \Bd D$. 
(This may be seen by an appeal to the maximum modulus
principle.) Below (except in \S\ref{quasipositive}, 
Remark \ref{Hartogs-Rosenthal}) 
we will assume $f_0(z)$ is a (non-zero) constant, 
that is, that there are no poles. 
This is only for convenience; everything
would work as well just assuming that no poles lie in $D$.

In \S\ref{braids} we recall the definition of \bydef{positive} 
closed braids, and define a strictly larger class, 
the \bydef{quasipositive} closed braids. The definition 
is purely braid-theoretic. Several mathematicians 
(including Murasugi, Stallings \cite{Stallings}, and Birman \cite{Birman}) 
have observed that many positive closed braids, 
in particular all those which are knots (rather than
links), are \bydef{fibred links}; 
there are quasipositive closed braids which are knots and not
fibred.

In \S\ref{quasipositive} we give one proof that the closed braid $K_\gamma$ 
is quasipositive. The proof is real semi-algebraic geometry, 
and gives a method (which is, alas, far from practicable in
most cases) of explicitly calculating the braid type of $K_\gamma$ 
in terms of one's knowledge of $\gamma$ and $f(z,w)$.

In \S\ref{analytic loops} we briefly discuss those loops in $M-V$, 
where $M$ is a simply-connected algebraic variety 
and $V$ is an algebraic subset, which are freely 
homotopic to loops which bound analytic (possibly 
singular) disks in all of $M$. In many cases, the free
homotopy classes of ``analytic boundaries'' turn out to be 
precisely those classes which are ``quasipositive'' in 
an appropriate sense. When $M$ is the space of unordered
$n$-tuples of (not necessarily distinct) complex numbers, 
and $V$ is the so-called ``discriminant locus'' of $n$-tuples 
with not all members distinct, the theory applies (to check 
one hypothesis, I use the method of 
\S\ref{quasipositive}), and we have the following 
theorem.

\begin{theorem*}
The closed braids $K_\gamma$ that arise from algebraic functions 
without poles are precisely the quasipositive closed braids.
\end{theorem*}

Here are some consequences of the theorem. 
Many more fibred links occur as $K_\gamma$ 
than just those associated to singular points 
of curves (as in \cite{Milnor})---these ``links of singularities'' 
may be recovered as a special case ($\gamma$ is a small 
circle enclosing a single point of $B$, for suitable $f(z,w)$). 
Many non-fibred knots and links occur as $K_\gamma$'s. 
And in each concordance class of links that appears at all, 
infinitely many distinct links occur; for instance 
(even for $f(z,w)$ as special as $w^3-3w+2z^m$, $m=1,2,3\dots$),
infinitely many distinct slice knots occur---a marked contrast 
to the links of singularities.

Remarks and examples conclude the paper.

\section{POSITIVE AND QUASIPOSITIVE BRAIDS AND CLOSED BRAIDS}\label{braids}

A general reference for the braid theory used here is \cite{Birman} 
(where a polyhedral approach is taken).

For $n\ge 2$ the algebraic $n$-string braid group $B_n$ 
is generated by $n-1$ standard generators $\s{1},\dots,\s{n-1}$
subject to the relations $\s{i}\s{i+1}\s{i}=\s{i+1}\s{i}\s{i+1}$
($i=1,\dots, n-2$), $\s{i}\s{j}=\s{j}\s{i}$ if $|i-j|>1$. 
A word $\s{k(1)}^{e(1)}\dotsm\s{k(m)}^{e(m)}$
(each $\e(j) = \pm 1$) in the generators and their inverses 
is \bydef{positive} if each $\e(j) = +1$, 
\bydef{strictly positive} if also every index from $1$ to $n-1$ 
occurs as some $k(j)$; an element $p$ of $B_n$ is (strictly) 
positive if it can be represented as a (strictly) positive word.

Let $K \sub \gamma \times \C$ be a closed braid in an open 
solid torus, with $K$, the simple closed curve $\gamma$, 
and $\C$ all oriented, and the projection from $K$ to $\gamma$ 
smooth and orientation preserving of degree $n$. 
It is well-known that the isotopy classes of such $K$ 
(say, ambient isotopy preserving the product structure 
of the solid torus) are in $1$-$1$ correspondence with 
conjugacy classes in $B_n$. The correspondence is implemented by
the choice of a diffeomorphism (preserving orientations) 
$h: \gamma \times \C \to S^1\times\R\times \R$ of the 
form $h(z,w) = (h_0(z), h_1(z,w), h_2(z,w))$ 
together with a basepoint $\exp\thinspace i\theta_0$ on $S^1$. 
Any such $h$ can be changed by an arbitrarily small isotopy, 
if necessary, to make it yield a ``good'' braid diagram 
$d(K)$ in the half-open rectangle 
$[\theta_0, \theta_0 + 2\pi] \times \R$ (project onto $S^1$
and take logarithms for the first coordinate, project onto 
the first $\R$ factor for the second coordinate, 
and at multiple points use the second $\R$ factor to determine 
under- and over-crossings)---``good'' in the sense that: 
$d(K)$ is the union of $n$ properly embedded arcs, 
on each of which the projection to 
$[\theta_0, \theta_0 + 2\pi]$ is a diffeomorphism;
there are no triple points of $d(K)$; 
there are only finitely many double points, all
interior to the rectangle, and at each of which 
the tangent lines to the two arcs are distinct; 
and the $\theta$ coordinates of distinct double 
points are distinct. From such a good
braid diagram $d(K)$ a word in the letters $\s{j}$ 
and their inverses may be read off, as follows. 
Let the $\theta$ coordinates of the double points 
be $\theta_1 < \theta_2 < \dotsm < \theta_m$. For each
$j =1, \dots, m$, there are precisely $n-1$ points in 
$\{\theta_j\} \times \R \cap d(K)$. Let the double point
be the $k(j)$th among them, in increasing order of $\R$ 
coordinate. Let $\psi_\theta$ and $\psi_2(\theta)$
parametrize the two arcs that cross at the double point 
in question, so labelled that $\psi'_1(\theta_j)>\psi'_2(\theta_j)$. 
Near $\theta_j$ there are smooth functions 
$\phi_1(\theta), \phi_2(\theta)$ so that 
$\theta \mapsto (\exp i\theta, \psi_l(\theta),\phi_l(\theta))$
($l = 1,2$) parametrize intervals on $h(K)$. 
Let $\e(j) =\sgn (\phi_2(\theta_j)-\phi_1(\theta_j)$. 
Then the word to be read off from $d(K)$ is 
$\prod\limits_{j=1}^m \s{k(j)}^{\e(j)}$.

A closed braid is \bydef{positive} if its corresponding 
conjugacy class in $B_n$ contains a positive braid. 
If $K$ has a braid diagram $d(K)$, as above, in which 
each exponent $\e(j)$ is $1$, certainly $K$ is positive.

Let $w_1,\dots , w_m$ be arbitrary words in 
$\s{1},\dots , \s{n-1}, \s{1}^{-1},\dots ,\s{n-1}^{-1}$. 
We will say that the word 
$w_1\s{k(1)}w_1^{-1}w_2\s{k(2)}w_2^{-1}\dotsm w_m\s{k(m)}w_m^{-1}$ is 
\bydef{quasipositive}, and that $q\in B_n$ is
\bydef{quasipositive} if it can be represented 
as a quasipositive word.

A closed braid is \bydef{quasipositive} if the corresponding 
conjugacy class in $B_n$ contains a quasipositive braid.

Now let $\gamma\times\C$ be embedded in $S^3$ 
as a tubular neighborhood of an unknotted circle,
and let $K$ be a closed $n$-string braid in that 
neighborhood. Corresponding to any good braid diagram 
$d(K)$, in which there are $m$ double points, 
there is a natural Seifert surface $S \sub S^3$ for K 
(i.e. an oriented surface with $\Bd S = K$) 
made up of $n$ disks connected by $m$ bands---the disks 
are ``stacked'' (they may be taken to be meridional
disks of the complementary solid torus to $\gamma \times\C$) 
and each band connects two adjacent disks in the stack, 
with a half-twist in one sense or the other depending on the sign
$\e(j)$ of the corresponding double point. 
(This construction by ``bands", following Murasugi, 
is expounded in Stallings's paper \cite{Stallings}. A 
general ``band representation'' which constructs 
``Seifert ribbons'' instead of Seifert surfaces, 
is discussed in \cite{Rudolph1981}.) 
As in \cite{Stallings}, when K is positive 
\textit{and so displayed by $d(K)$}, any connected component 
$S_0$ of $S$ has the property that the push-off map 
$\pi_1(S_0)\to \pi_1(S^3-S_0)$ (defined by taking a
nowhere-zero normal vectorfield on $S_0$ and using it 
to push any loop on $S_0$ into the complement of $S_0$) 
is a bijection. It then follows from a theorem of Neuwirth and
Stallings that the boundary of $S_0$, a union of components 
of the link $K$, \textit{is a fibred link}.
In particular, $K$ is fibred if either $K$ is a knot 
or $S$ is connected, which last happens if and only 
if the word of $d(K)$ is strictly positive. 
Details of the proof appear in \cite{Birman-Williams}.
\appnote{fibration of positive braid}

\section{THE CLOSED BRAIDS \protect$\mathbf{K_\gamma}$ ARE QUASIPOSITIVE}
\label{quasipositive}

Until further notice, our algebraic functions will not have any poles.

Let $\pi = \pr_1|V_f: V_f\to \C$. We begin by observing that there is no 
loss of generality, for the purposes of studying all the braids $K_\gamma$, 
in assuming that $V_f$ is a non-singular
curve and that for each $z_0 \in B$, the fibre 
$\pi^{-1}(z_0)$ consists of $n-1$ distinct points, at
one of which $V_f$ has a vertical tangent. 
Indeed, if this is not so already, any sufficiently 
small change in the constant term of $f_{n-1}(z)$ will 
make it so; while the closed braids lying over a fixed $\gamma$  
on the two curves $V_f$  and $V_{f+\e w}$ are surely isotopic
(by a vertical isotopy) for all sufficiently small $\e$.

Now suppose that $\gamma_0$ and $\gamma_1$ are isotopic 
in the complement of $B$. The differential
$D\pi$ is surjective off $\pi^{-1}(B)$; 
so the isotopy lifts to an isotopy of embeddings between
$K_{\gamma_0}\hookrightarrow \gamma_0 \times \C$
and $K_{\gamma_1} \times \C$. In the special case 
that $\gamma_0$ and $-\gamma_1$ cobound an annulus
$A$ in the complement of $B$, then the union of 
annuli $\pi^{-1}(A) \sub V_f$  is the trace of an
isotopy between the closed braids.

To show that $K_\gamma$ is quasipositive we will 
isotope $\gamma$ to a more-or-less normal form
for which the conclusion will be obvious. 
(All the unbridgeable gap between positive
and quasipositive lies in that ``more-or-less''!)

We will begin by constructing an oriented graph 
(smoothly embedded in the plane) with vertices 
including all the points of $B$.
\appnote{Orevkov} %
Let $z_1, \dots, z_l$
be the points of $B$, and for $j = 1, \dots , l$ let 
$w_{j,l}, \dots , w_{j,n-1}$ be the $n-1$ distinct roots 
of $f(z_j,w) = 0$.  Then for all but
finitely many $\theta\in [0, 2\pi]$ the $n-1$ real 
numbers $\Re((\exp i\theta)w_{j,k})$, 
$k = 1,\dots , n-1$, are pairwise distinct, 
for each $j = 1, \dots, l$. 
Changing the $w$-coordinate by a rotation, then,
we may assume without loss of generality that $\theta = 0$ works, 
that is, that at each point $z_j$ the $n-1$ real parts $\Re w_{j,k}$ 
are pairwise distinct. Let 
$B^{+} = B \cup \{z\in \C-B\Suchthat
\text{ for some two distinct solutions } w_1, w_2 
\text{ of } f(z,w) = 0, \Re w_1 = \Re w_2\}$. 
Then $B^{+}$ is the projection of a real algebraic set, 
so on general principles it is a real semialgebraic set,
evidently of dimension $1$, and so a graph; we will see 
this directly in the course of establishing its local 
structure. We will find a locally-finite (actually finite) 
subset $B_0$ of $B^{+}$, containing $B$, 
so that $\C$ is stratified by $B_0$, $B^{+}-B_0$, 
$\C-B^{+}$. Let us consider the
intersection of $B^{+}$ with a disk around an arbitrary point 
of $\C$. If this point $z_0$  does not
belong to $B$, let $\e > 0$ be sufficiently small 
that the disk $D_\e(z_0)$ is disjoint from $B$. Then
on this disk there are analytic functions 
$w_j(z)$ so that $\pi^{-1}(D_\e(z_0))$ is the union of the
graphs of the functions $w_j$.  Thus $B^{+} \cap D_\e$
is the union of sets $A_{j,k} = \{z\in D_\e(z_0)\Suchthat 
\Re(w_j(z)-w_k(z)) = 0\}$. Each difference $w_j-w_k$ is analytic, 
not identically $0$, and so, near any point of $D_\e(z_0)$, 
$w_j - w_k$ is a branched cover of its image; so the real analytic
set $A_{j,k}$ is a $1$-complex, smoothly embedded near its 
manifold points, and near its finitely many non-manifold points 
(which we assign to $B_0$) smoothly equivalent to a union of 
diameters in a disk. Likewise, distinct sets $A_{j,k}, A_{g,h}$ cross 
only finitely often; put their intersections in $B_0$ too.

If we look near a point $z_j$ of $B$ the situation is slightly different. 
Here, for small $\e > 0$, $\pi^{-1}(D_\e(z_j)$ consists of not 
$n$ but $n-1$ smooth disks. There are $n-2$ functions
$w_k(z)$ analytic on $D_\e(z_j)$ whose graphs are 
$n-2$ of these disks; the last disk is
parametrized by $t\mapsto (z_j+t^2, w(t))$, 
where $|t|^2 < \e$, $w(t)$ is analytic, and 
$w'(0)\ne 0$ (we are at a simple vertical tangent). 
Since we have assumed $\Re w_i(z_j),\dots , \Re w_{n-2}(z_j), \Re w(0)$ 
are distinct, after possibly shrinking $\e$ we can guarantee 
that $B^{+}\cap D_\e(z_j)$ has no contributions from the 
interaction of any of the $w_k(z)$ with each other or with $w(t)$:
we will have simply $B^{+} \cap D_\e(z)=\{z_j+t^2\Suchthat %
|t|^2<\e, \Re(w(t)-w(-t)=0\}$. But, like
$w(t)$, $w(t)-w(-t)$ has non-zero derivative 
at $t=0$, so (shrinking again if necessary)
we see that $\{t\Suchthat |t|^2 < \e, \Re(w(t)-w(-t) = 0\}$ 
is smoothly (and equivariantly) equivalent to a diameter of 
the $t$-disk, and its image in $B^{+}$ is smoothly equivalent 
to a radius of $D_\e(z_j)$.

We now orient $B^{+}$, at the same time labelling each edge 
with one of the symbols $\s{1},\dots, \s{n-1}$. Let $A$ be an arc 
in $B^{+}-B_0$. Then anywhere in the interior of $A$, one may
find a short transverse arc which intersects $A$ only in one point, 
and $B^{+}$ nowhere else. Over such an arc the $n$ branches of $w(z)$ 
are distinct, and even their real parts are distinct except where the 
transverse arc crosses $A$: at that point, for some $k$,
$1\le k\le n-1$, the branches with real parts 
$k$th-greatest and $(k+1)$st-greatest among all
the branches have equal real part; label $A$ with 
$\s{k}$. (Clearly this label is independent of
the transverse arc.) Orient $A$ so that, when the 
orientation of the transverse arc, following the 
orientation of $A$, gives the complex orientation of $\C$, 
the braid diagram over the transverse arc is one for 
$\s{k}$ (rather than for $\s{k}^{-1}$).

Let $\gamma$ be a smooth simple closed curve in 
$\C-B$, oriented counterclockwise, and bounding the 
bounded region $D$.  Let $z_1,\dots, z_s$ be the points 
of $B \cap D$, let $D_j = D_\e(z_j)$,
and let $C_j =\Bd D_j$ oriented counterclockwise, 
for $j = 1,\dots, s$. For sufficiently small $\e$ the disks 
$D_j$ lie in $D$ and are pairwise disjoint. By 
a traditional construction of the theory of algebraic functions, 
there is a disk $D_0 = D_{\e_0}(z_0) \sub %
D - \bigcup\limits_{j=1}^s D_j$ with boundary
$C_0$ (oriented counterclockwise), and pairwise 
disjoint smooth embeddings $a_j: [0,1]\to D$
$(j = 1, \dots, s)$ with $a_j(0)\in C_0, a_j(1) 
\in C_j, a(]0,1[) \sub D - \bigcup\limits_{k=0}^s D_k$, and 
$a_j$ perpendicular to $C_0$ and $C_j$ at its ends, all 
so that $\gamma$ is isotopic in $D - B$ to a simple closed
curve $\gamma'$ which ``follows the arcs and circles.'' 
Formally, $\gamma' = \Bd (D_0 \cup \bigcup\limits_{j=1}^s N_j %
\cup \bigcup\limits_{j=1}^s D_j)$, where the sets $N_j$ are 
``strips"---pairwise disjoint product neighborhoods of the arcs
$a_j([0,1])$, say $N_j = \nu_j([-1,1] \times [0,1])$, 
where $\nu_j$ is an embedding such that $\nu_j(0,t)=a_j(t)$
$(t \in [0, 1])$, etc.

Now we involve $B^{+}$. Without loss of generality, 
we assume that $D_j$ ($j = 1,\dots , s$) intersects 
$B^{+}$ only in an arc that joins $z_j$ to $C_j$, 
and that $D_0$ is disjoint from $B^{+}$. 
It is clear that, in performing the traditional construction, 
we may so arrange things that the embeddings $a_j$ 
are transverse to the stratification---they miss $B_0$
and cross the manifold points of $B^{+}$ transversely in 
the ordinary sense---and then make the product neighborhoods 
$N_j$ so narrow that $N_j \cap B^{+}$ is itself a product 
$[-1,1]\times (a_j([0,1]) \cap B^{+})$.

Let $h_0: \gamma'\to S^1$ be a diffeomorphism so that 
$h_0^{-1}(1)$ is a point on $C_0$; define 
$h: \gamma' \times \C \to  S^1 \times \R\times \R$
by $h(z,w) = (h_0(z), \Re w, \Im w)$. 
I claim that applying the construction of \S\ref{braids}
to this $h$ (with base-point $1$ on $S^1$) yields a 
good braid diagram $d(K_\gamma')$ for which the braid 
word is already in the form 
$\prod\limits_{j=1}^m \a_j\s{k(j)}\a_j^{-1}$; 
so that $K_\gamma'$ and $K_\gamma$ are
quasipositive. Indeed, the diagram $d(K_\gamma')$ 
is the ``product'' in an obvious sense of diagrams 
for the (non-closed) braids which correspond to the 
successive arcs $\nu_{j(1)}(\{1\} \times [0,1])$, 
$C_{j(1)}-\nu_{j(1)}({]}-1,1{[}\times\{1\})$, 
$\nu_{j(1)}(\{-1\} \times [0,1]), \dots$  of $\gamma'$ 
(where the order in which the points of $B\cap D$ are gone 
around is $z_{j(1)},\dots, z_{j(s)}$, and where the arcs
$\nu_{j(k)}(\{-1\}\times [0,1])$ are of course traversed 
from the $1$ end to the $0$ end). Each arc
contributes, in turn, the word in the symbols 
$\s{1},\dots, \s{n-1}, \s{1}^{-1},\dots, \s{n-1}^{-1}$
which is given by its successive crossings of the labelled 
arcs of $B^{+}-B_0$ (a crossing which, following the 
orientation of the arc, gives the wrong orientation to 
$\C$, is what merits the exponent $-1$). 
Obviously, by our construction, the two edges of a strip 
$N_j$ give (up to orientation) the same word as the central 
arc $a_j([0,1])$, call it $\a_j$.  So the claim of
quasipositivity is proved once one sees that the diagrams 
corresponding to the arcs on the circles $C_j$ ($j = 1, \dots , s$) 
contribute exactly a generator $\s{k(j)}$, and not the inverse of
a generator. (Certainly by construction each such arc meets $B^{+}$ 
in just one point.) The exponent is seen to be $+1$ in all cases; 
it suffices to study just one example, for instance 
$f(z,w) = w^2 + z$, where $B=\{0\}$, $B^{+}$ is the 
non-negative real numbers, and the conclusion is obvious.

We have proved that \textit{if $f(z,w)$ has no poles 
inside $\gamma$, the closed braid $K_\gamma$ is quasipositive.} 
A converse will be proved in the next section.

\begin{remarks}
\begin{inparaenum}
\item\label{exponent sum}
The \bydef{exponent sum} $e(w)$ of a braid word 
$\mathop{\scriptscriptstyle\mathop{\textstyle\Pi}\limits_{j=1}^m}%
\s{k(j)}^{\e(j)}$ is 
${\scriptscriptstyle\mathop{\textstyle\Sigma}\limits_{j=1}^m} \e(j)$.
From the form of the relations in $B_n$, this is 
actually defined on braids; clearly it is
conjugation invariant, so it is an isotopy invariant 
of closed braids. The exponent sum of a quasipositive 
braid is non-negative. The proof above actually shows that the
exponent sum of $K\gamma$ is the number of points of $B$ enclosed 
by $\gamma$ (counting multiplicities appropriately if $f(z,w)$ 
is not restricted to simple vertical tangents and no singularities). 
It is easy to see that the exponent sum of a closed braid $K$ equals 
$sw(K)$, the \bydef{self-winding} defined by Laufer \cite{Laufer}. 
The proof above readily generalizes to analytic (rather than 
simply algebraic curves), and some theorems of \cite{Laufer} can be recovered
quickly.

\item\label{Hartogs-Rosenthal}
We have excluded from consideration simple closed curves 
$\gamma$ enclosing poles of our algebraic function. This 
is because, on the one hand, if $\gamma$ does enclose any
poles of $f(z,w)$ then the closed braid $K_\gamma$ is not 
the whole intersection of $V_f$ with $\Bd (D \times D_r^2)$
for any $r$---there are always components in $D \times \Bd D_r^2$ 
corresponding to the poles; while, on the other hand, if 
we allow poles then \textit{every isotopy class of closed braid}
can be realized as the braided part $K_\gamma$ of that intersection, 
for appropriate $f(z,w)$ and $\gamma$.  The proof is by the theory 
of rational approximation. Let $\gamma = \{z\Suchthat |z| = 1\}$. 
Let $K_0 \sub \gamma \times \C$ be a closed braid, 
not necessarily smoothly embedded, with components 
$C_1,\dots, C_d$ of degrees $n_1,\dots , n_d$. 
For a suitable large constant $M$, the polynomial 
$p(z) = M(z-1)^{n_1} \dotsm (z-d)^{n_d}$ is such that 
the compact set $P = \{z: |p(z)| \le 1\}$ is the union 
of $d$ components, each diffeomorphic to a disk, on the
boundaries of which $p(z)$ has degrees $n_1,\dots, n_d$ 
respectively. Evidently, there is a unique continuous function 
$q_0(z)$ defined on $\Bd P$ such that the pair 
$(p, q_0): \Bd P \to \gamma \times \C$ parametrizes 
$K_0$. According to the Hartogs-Rosenthal Theorem \cite{Gamelin}, 
on any compact subset of $\C$ with measure $0$ (e.g. $\Bd P$) 
the rational functions with poles off the compact set are 
uniformly dense in the continuous functions. Let $q(z)$ be a 
rational approximation to $q_0(z)$ so close that 
$K = (p,q)(\Bd P)$ lies inside a tubular neighborhood of
$K_0$ in $\gamma \times \C$ (which exists, even 
though $K_0$ may not be smooth, because $K_0$ is a closed
braid); then $K$ and $K_0$ are isotopic (by a vertical isotopy). 
But $(p,q)(\C) = V$ is an algebraic curve in $\C^2$ 
(generally with many singularities), that is, $V = V_f$  for some
$f(z,w)$.

Of course, when $q$ has poles interior to $P$ as well 
as in $\C-P$, there will be poles of $f(z,w)$ enclosed by $\gamma$.

\item\label{later use}
For later use, and intrinsic interest, we give some calculations 
of sets $B^{+}$ in particular examples.
\end{inparaenum}

\begin{example}
$f(z,w) = w^2-z$. Here $w_1 = \sqrt{\vphantom{(}}(z)$, 
$w_2 =-\sqrt{\vphantom{(}}(z)$, and $\Re w_1 = \Re w_2$ iff
$w_1$ and $w_2$ are pure imaginary iff $z$ is negative real; 
thus $B^{+}$ is the ray ${]}-\infty, 0{]}$ ending in
$0$, the only point of $B$; the ray is oriented away from $0$, 
and labelled $\s1$. More
generally, if $f(z,w) = w^2-z^n-1$, then 
$B^{+} = \{z\Suchthat z^n + 1 \text{ is negative real}\}$ 
is the union of $n$ rays, oriented outward, emanating 
from the $n$th roots of $1$, all labelled $\s1$. Of course,
in the $2$-string braid group, which is infinite cyclic, 
quasipositive is the same as positive.
\end{example}

\begin{example}
$f(z,w) = w^3-3w + 2z^n$. 
If $w_1$, $w_2$, and $w_3$ are the three roots of
$f(z,w)=0$, then $w_1+w_2+w_3=0$, $w_1 w_2+w_1 w_3+w_2 w_3= -3$, 
and $w_1 w_2 w_3= -2z^n$.  Eliminating $w_3$ between the first 
two equations, we get the quadratic relation 
$w_2^2 + w_1 w_2 + (w_1^2-3) = 0$, whence $\{w_2, w_3\}=%
\{\frac12(-w_1 + \sqrt{\vphantom{(}}(-3w_1^2+ 12),\frac12(-w_1-\sqrt{\vphantom{(}}(-3 w_12 + 12\}$.
The indices are irrelevant; there is perfect symmetry, and
we see that $B^{+}=\{z\Suchthat \Re w_2= \Re w_3\}=%
\{z\Suchthat \sqrt{\vphantom{(}}(-3w_1^2+12) \text{ is pure imaginary}\}=%
\{z\Suchthat -3w_1^2 \in {]}-\infty, -12{]}\}$. 
For $n=1$, $B^{+}$ is thus the two rays ${]}-\infty,-1]$ and 
$[1,\infty{[}$; in general, $B^{+}$ is the union of $2n$ rays, 
oriented outward, emanating from the $2n$th roots of $1$, 
and labelled alternately $\s2$ and $\s1$. For $n=4$, we get 
an example of a quasipositive, not positive, knot $K_\gamma$ 
for the curve pictured in Figure~\ref{Figure 1}; the braid 
word here is $\s1\s2^3\s1\s2^{-3}$. This knot is $8_{20}$ of the 
Alexander-Briggs table; it is slice---indeed, ribbon---and 
non-trivial; it cannot be positive because, for instance, 
according to \cite{Rudolph1982} a non-trivial positive closed 
braid has signature greater than $0$.
\end{example}

\begin{figure}
\centering
\includegraphics[width=2.0in]{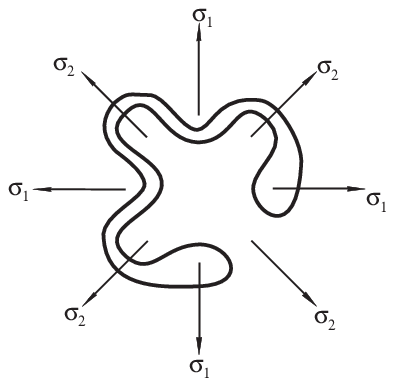}
\caption{\label{Figure 1}}
\end{figure}

\begin{example}\label{Example 3.3}
(This example will be used in the next section to establish that all
quasipositive closed braids occur as $K_\gamma$'s.)  Consider the 
reducible polynomial $f(z,w)=P(w)(w-z)$, where $P$ is a polynomial 
in $w$ without double roots. Here $B={z\Suchthat P(z)=0}$ is just 
the set of roots of $P$, and $B^{+}$ will either be all of $\C$ 
(in the unfortunate case, ruled out in the discussion above by a 
rotation of $w$ when necessary, that some two distinct roots of $P$ 
have equal real part) or, generically, the union of $n$ straight 
(real) lines $\Re z=r_j$ ($j=1,\dots, n$), where $r_j$ is the real 
part of a (unique) root of $P$: $B_0$ here is just $B$. Now suppose 
$P$ has real coefficients, and consider, for $\e\ne0$ small and real, 
the set $B_\e^{+}$ corresponding to $f(z,w)+\e$, and its
distinguished subset $B_\e$.  Evidently these sets are invariant 
under complex conjugation of the variable $z$. One sees that, 
in fact, the points of the original $B$ were ``to be counted twice'' 
and that as $\e$ moves away from $0$ these points of multiplicity two
alternately (with increasing $r_j$) bifurcate to two real points and 
to two conjugate, non-real points. Further, it is not much harder 
to see that the interval of the real line between the points of a 
real pair itself lies entirely in $B_\e^{+}$. Only in the simplest 
case, when $P$ is linear, have I been able to get an explicit 
description of the full set $B_\e^{+}$; but this suffices to give 
an adequate qualitative description in the general case. Namely, if
$P(w)=w$, say, then $B_\e^{+}=\{z\Suchthat w^2-wz + \e\} \text{ has 
two real roots with equal real parts}\} =%
\{z\Suchthat \sqrt{\vphantom{(}}(z^2-4\e) \text{ is pure imaginary}\}=%
\{z\Suchthat z^2\in {]}-\infty,4\e]\}$. When $\e<0$, 
this is the union of two rays lying on the imaginary axis, 
oriented outward; when $\e>0$, however, it is a
cross, containing the whole imaginary axis and a short 
interval of the real axis---the short arms oriented towards 
the crossing point, the long arms out to infinity. Now for
a polynomial $P$ of higher degree, there is a neighborhood 
$N$ of $B$ which is a union of disjoint disks around the 
roots of $P$, so that for $\e$ sufficiently small (and real) 
the set $B_\e^{+}$ looks like $B^{+}$ outside $N$ (that is, 
it consists of two proper arcs leaving each disk of $N$ and 
going to infinity without crossing) while inside alternate 
disks of $N$ (from left to right) $B_\e^{+}$ looks like the 
case $P(w)=w$, with an $\e$ of the same or opposite sign. 
So the whole set $B_\e^{+}$ is, qualitatively, a sequence 
of alternate crosses and double-rays; Figure~\ref{Figure 2}
gives a sketch in case $P(w)=w(w-1)(w+1)$. 
The orientations are as in the linear
model, and from left to right the arcs of $B_\e^{+}$ are labelled 
$\s1,\dots,\s{n}$ (where $n$ is the degree of $P$) in batches. 
For later use note that, from an arbitrary basepoint $*$ (off $B^{+}$)
for each $j$ a loop can be drawn whose word in the labels 
$\s{j}$ and $\s{j}^{-1}$ is freely equal (in the free group 
on the labels) to $\s{j}$. For instance, for $*$ to the far 
left in Fig.~2, a loop for $\s1$ is obvious; a loop for 
$\s2$ can slip between the two rays labelled $\s1$, do the
obvious, and slip back; a loop for $\s3$ will have to intersect 
the cross labelled $\s2$, but if it goes through the gap 
between the two ends of the short arm it will pick up
successively $\s2$ and $\s2^{-1}$; and so on.
\end{example}

\begin{figure}
\centering
\includegraphics[width=2.0in]{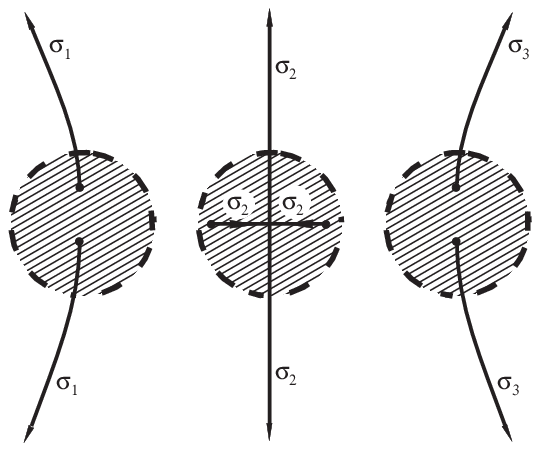}
\caption{\label{Figure 2}}
\end{figure}
\end{remarks}

\section{ANALYTIC LOOPS IN THE CONFIGURATION SPACE}\label{analytic loops}
Throughout this section let $D$ be the closed unit disk in $\C$, 
$S^1=\Bd D$ its boundary oriented counterclockwise.

If $X$ is a complex analytic space, an 
\bydef{analytic disk} in $X$ is a map $i: D \to X$ which is
the restriction to $D$ of a complex analytic map on some 
slightly larger open disk; an \bydef{analytic loop} is the 
oriented boundary of an analytic disk. Suppose $X$ is simply
connected, and $V \sub X$ is a closed analytic subset such 
that $X-V$ is connected but no longer simply connected. We may 
ask, which non-trivial homotopy classes of loops in $X-V$ 
contain representatives which are analytic loops in $X$?

Even when the question is asked in such generality, partial 
answers can be given.  For our present purposes, however, it 
is enough to have the answer with $X$ and $V$
considerably restricted. So, let $X=\C^n$ be affine space, 
and let $V \sub \C^n$ be an algebraic hypersurface 
$V=V_f=\{\mathbf{z} \in \C^n\Suchthat f(\mathbf{z})=0\}$, 
possibly singular and/or reducible (but without multiple 
components). The complex manifold $R(V)$ of regular points of 
$V$ is of (real) codimension $2$ in $\C^n$, and is everywhere 
dense in $V$; let its connected components be $R_1,\dots, R_s$. 
For some arbitrary point on each $R_i$, let $D_i$ be an oriented
normal $2$-disk intersecting $V$ only at that point, 
and there positively (with respect to the complex orientations 
of $R(V)$ and $\C^n$); for some fixed basepoint $*$ not on $V$, 
let $a_i$ be an arc in $\C^n-V$ from $*$ to a point on 
$\Bd D_i$; let $l_i$ be a loop which runs from $*$ along $a_j$ 
to $\Bd D_i$, once around $\Bd D_i$ countercloskwise, and back 
along $a_j$ to $*$; and let $[l_i]$ be the class of $l_i$ in 
$\pi_1(\C^n-V)$; all for $i=1,\dots, s$. 
For later use, in the particular case that $n=1$ and $V$ is a 
finite set of points, each one a component $R_i$, let us demand 
further that the disks $D_i$ be pairwise disjoint from each other 
and from $*$, and that the arcs $a_i$ be simple, pairwise disjoint 
except for their common endpoint $*$, and outside the union of 
the $D_i$ (except for their other endpoints).

An element of $\pi_1(\C^n-V;*)$ which can be written as a product 
$\mathop{\scriptscriptstyle\mathop{\textstyle\Pi}\limits_{j=1}^m}%
w_i[l_j(i)]w_i^{-1}$ of conjugates of the classes $[l_i]$ will be 
called a \bydef{quasipositive} element of the fundamental
group. Quasipositivity is invariant under conjugation, and thus 
is really a property of free homotopy classes of loops.

\begin{lemma}\label{Lemma 1}
An analytic loop in $\C^n-V$ represents a quasipositive conjugacy 
class in $\pi_1(\C^n-V;*)$.
\end{lemma}

\begin{proof}
Let $i: D \to \C^n$ be an analytic disk in $\C^n$ 
with $i(S^1) \cap V=\emptyset$. Replacing $i$ by a
sufficiently close approximation (for instance, a high-order 
Taylor polynomial at $0$) we may assume $i$ is the restriction 
to $D$ of a (vector-valued) complex polynomial $p(t)$ of a
single complex variable t, without changing the (free) 
homotopy class of $i(S^1)$ in the complement of $V$. 
In $\C\times\C^n\times\C^n$ let $Z$ be the set 
$\{(t,\mathbf{\e},\mathbf{z})\Suchthat 
\mathbf{z}=p(t)+\mathbf{\e} \text{ belongs to } S(V)\}$,
where $S(V)=V-R(V)$ is the singular set of $V$, an algebraic 
set of complex dimension no greater than $n-2$. 
Then $Z$ is an algebraic subset of $\C^{2n+1}$. 
Its complex dimension is no greater than $n-1$, 
for $\mathbf{z}$ varies in a set of dimension at most $n-2$,
$p(t)$ is on a curve, and $\mathbf{\e}$ is determined by 
$\mathbf{z}$ $p(t)$. Then the projection of $Z$ onto
the second factor, $pr_2(Z) \sub \C^n$, is again an algebraic 
set of dimension at most $n-1$.
Then almost any $\e$, in particular, almost any $\e$ 
sufficiently close to $\mathbf{0}$, is not in $pr_2(Z)$.
Translating $i(D)$ by an appropriate small $\mathbf{\e}$ 
will not change the free homotopy class of the analytic 
loop $i(S^1)$ while ensuring that $p(\C)$ and its subset 
the new analytic disk meet $S(V)$ nowhere.  Now the whole 
intersection of the analytic disk and $V$ is in the
manifold $R(V)$ and it is a simple matter to make the 
intersection transverse, when it will appear that each 
point of intersection counts $+1$ because $p(\C)$ and 
$R(V)$ are complex manifolds. Since the boundaries of 
two normal disks (positively oriented) at any two points 
of a component $R_i$ are freely homotopic, the analytic 
loop is a product of conjugates of the loops $l_i$.
\end{proof}

\begin{lemma}\label{Lemma 2}
Conversely, when $n=1$, every quasipositive conjugacy class 
in $\pi_1(\C-\{z_1,\dots, z_s\})$ is represented by an analytic 
loop in $\C$.
\end{lemma}

I do not know if Lemma~\ref{Lemma 2} is true when 
$n\ne 1$. However, the following immediate consequence of 
Lemma~\ref{Lemma 2} suffices to replace the putative stronger 
version for our purposes.

\begin{unproved}{cor}
If there is a proper analytic map $L$ of $\C$ into $\C^n$ so that 
the induced homomorphism $\pi_1(\C-L^{-1}(V))\to \pi_1(\C^n-V)$ is 
surjective, then every quasipositive conjugacy class in $\pi_1(\C^n-V)$ 
is represented by an analytic loop (which in fact bounds an analytic 
disk lying on $L(\C))$.
\end{unproved}

\begin{proof}[Proof of Lemma~\protect\ref{Lemma 2}.] 
Let $\a=\mathop{\scriptscriptstyle\mathop{\textstyle\Pi}\limits_{j=1}^m}%
w_i[l_j(i)]w_i^{-1}\in \pi_1(\C=\{z_1,\dots,z_s\},*)$ be quasipositive.
Let the disks $D_j$ ($j=1,\dots, s$) be as above, 
let $D_0$ be a disk centered at $*$ and disjoint
from all the other $D_j$, and suppose for neatness 
that for each $j=1,\dots,s$ the arc $a_j$ intersects $D_0$ 
in a radius of $D_0$, and comes into $D_j$ normally. 
Let $c(j)\ge 0$ be the number of times the index $j$ appears as 
$j(i)$ in the given presentation of $\a$, as $i$ runs
from $1$ to $m$. Let $D'_{j,c}$ ($j=1,\dots, s, c=1,\dots, c(j)$) 
and $D'_0$ be $2$-disks which we think of as ($2$-dimensional) 
$0$-handles, and let $N_i$ ($i=1,\dots, m$) be strips, 
each homeomorphic to $[-1,1]\times [0,1]$, which we think of as 
$1$-handles. Fix orientations on all the handles. Take $m$ 
disjoint closed intervals, successive in the cyclic order, 
on $\Bd D_0$, and one closed interval on each of the 
$\Bd D'_{j,c}$ (of which there are $m$ all together). 
We form an identification space from the disjoint union 
of all the $0$- and $1$-handles as follows: orientedly, 
attach one end $[-1,1]\times \{0\}$ of $N_i$ to the $i$th 
chosen interval on $\Bd D_0$, and the other end 
$[-1,1] \times \{1\}$ to the chosen interval on 
$\Bd D'_{j(i),c}$ (where $c$ is the number of $k$ with 
$k\le i$, $j(k)=j(i)$). Then this identification space $D''$ 
is homeomorphic to a disk. We will map $D''$ into $\C$ 
handle by handle. First each $D_{j,c}$ is mapped homeomorphically, 
preserving orientation, onto $D_j$ so that the image of the chosen 
interval on $\Bd D_{j,c}$ is centered at the end of $a_j$ on $\Bd D_j$; 
and $D_0$ is mapped homeomorphically, preserving orientation, 
onto $D_i$. For each conjugator $w_i$, find an immersed arc 
in $\C$ which begins (outward normal) in the image on 
$\Bd D_0$ of the $i$th chosen interval on $\Bd D_0$ 
and represents $w_i$ in $\pi_1(\C-\{z_1,\dots, z_s\}, D_0)$; 
then map the center line $\{0\}\times [0,1]$ of $N_i$ to an arc 
which follows the arc representing $w_i$ from $\Bd D_0$ 
back to $D_0$, then in $D_0$ to $*$, and then along
$a_{j(i)}$ to $D_j$.  Because the exponent of $[l_{j(i)}]$
in $\a$ is $+1$ and not $-1$, the map on this center line can 
be extended over all of $N_i$ to give an immersed tubular 
neighborhood of the image of the centerline, which respects the 
identifications at both ends.  The map so constructed is an 
immersion on the interior $\Interior{D}''$, and on the boundary
represents $\a$. By ``transport of structure'' the interior of 
$D''$ becomes a Riemann surface, and by the Riemann Mapping theorem 
there is an analytic homeomorphism 
${\Interior{D}}_{1+\e} \to {\Interior{D}}''$, 
where ${\Interior{D}}_{1+\e}=\{z\Suchthat |z|<1+\e\}$, 
for any $\e>0$. For appropriately small $\e$, if $i$ is the composite 
$D \sub \Interior{D}_{1+\e} \to \Interior{D}'' \to \C$, 
then $i$ is an analytic disk whose boundary $i(S^1)$
represents (the conjugacy class of) $\a$. (A tiny bit more 
juggling could assure that $i(S^1)$ passed through $*$.)
\appnote{Blank}
\end{proof}

Presumably the hypothesis of the corollary is always true, 
even with $L$ a linearly parametrized straight line in 
sufficiently general position (see \cite[p.~33]{Lefschetz}). 
In any case, consider the following example.

\begin{example}
The group $B_n$ may be defined topologically as the fundamental group
of the \bydef{configuration space} of unordered $n$-tuples of 
distinct points in $\R^2$. Reading $\C$ for $\R^2$, one may recognize 
that, first, the space $\C^n/\mathscr{S}_n$ (where $\mathscr{S}_n$,
the symmetric group on $n$ letters, acts by permuting the coordinates) 
of unordered $n$-tuples of complex numbers (distinct or not) 
is in a natural way equal to $\C^n$ again, by the theorem on
symmetric polynomials; and, second, that the so-called 
``multi-diagonal'' or discriminant locus, consisting of 
unordered $n$-tuples of which two (at least) are equal, 
is an algebraic hypersurface $V_\Delta$ in the affine space 
$\C^n/\mathscr{S}_n$.  I claim that Example~\ref{Example 3.3} 
provides one with a line $L$ in $\C^n/\mathscr{S}_n$ satisfying 
the hypothesis of the Corollary to Lemma~\ref{Lemma 2}. 
For, what ``is'' an element of $\C^n/\mathscr{S}_n$ 
but the monic polynomial of degree $n$, in one complex variable $w$, 
whose roots are the unordered $n$-tuple in question? 
Under this identification, the affine coordinates in $\C^n/\mathscr{S}_n$ 
are precisely the significant coefficients of that polynomial 
(to wit, up to sign, the elementary symmetric functions of the 
roots). Now, if the polynomial $P(w)$ in Example~\ref{Example 3.3}
is chosen monic of degree $n-1$, 
then the assignment $L: z\mapsto P(w)(w-z) + \e \in C[w]$
of a monic polynomial of degree $n$ is clearly a linear 
parametrization of a straight line in $\C^n/\mathscr{S}_n$.
The work done in the example shows that 
$\pi_1(L(\C)-V_\Delta) \to \pi_1(\C^n/\mathscr{S}_n-V_\Delta)=B_n$
is surjective. Further, the two uses of the word ``quasipositive'' 
coincide here.
\end{example}

According to this example and the corollary, every quasipositive 
element of $B_n$, when considered as a homotopy class in the 
configuration space, contains an analytic loop in 
$\C^n/\mathscr{S}_n$.  But an analytic disk $i: D\to \C^n/\mathscr{S}_n$
is nothing more nor less than an $n$-valued analytic function on $D$, 
that is, an analytic subset of $D \times \C$ which projects
properly and $n$-to-$1$ (counting multiplicities) to $D.$ 
Without changing the free homotopy class of $i(S^1)$ in 
$\C^n/\mathscr{S}_n-V_\Delta$, one may (as in the proof 
of Lemma~\ref{Lemma 1}) replace the analytic function by 
(the restriction to $D$ of) a vector-valued polynomial; and a
polynomial map from $\C$ to $\C^n/\mathscr{S}_n$ is precisely 
an $n$-valued \textit{algebraic} function without poles. We 
have proved the following.

\begin{theorem*}
The closed braids that arise from algebraic functions without poles are
precisely the quasipositive closed braids. 
\end{theorem*}

\begin{remarks}
\begin{inparaenum}
\item
Which classes in $\pi_1(X-V;*)$ are represented by analytic loops
depends not only on $X-V$ but very strongly on $X$ as well. 
For instance, the natural way to complete the affine space 
$\C^n/\mathscr{S}_n$ is to $(\C\P^1)^n/\mathscr{S}_n$, which 
is canonically $\C\P^n$. Let $\bar{V}_\Delta$ be the completion
of $V_\Delta$ in $\C\P^n$ and let ${\C\P_\infty}^{n-1}$ be 
$\C\P^n-\C^n$, that is, the unordered $n$-tuples of extended 
complex numbers one at least of which is $\infty$. Then certainly 
$(\C^n/\mathscr{S}_n)-V_\Delta = 
((\C\P^1)^n/\mathscr{S}_n)-(\bar{V}_\Delta\cup {\C\P_\infty}^{n-1})$.
But the loops in this space, which are boundaries of analytic disks
in the whole projective space, fall into every homotopy class:
everything is quasipositive.  Indeed, an analytic disk in 
the projective space is an $n$-valued analytic function 
\textit{with poles allowed}; the poles correspond to intersections
of the disk with ${\C\P_\infty}^{n-1}$. Then by 
Remark~\ref{Hartogs-Rosenthal} of \S\ref{quasipositive} 
we actually have that any loop at all can be perturbed by an arbitrarily 
small amount, to become the boundary of an analytic disk 
(probably crossing infinity). In general, it appears that 
there will be more analytic disks in a projective variety 
than in a comparable affine one.
\item
If $X$ is a simply connected complex manifold, 
and $V$ is a non-singular analytic subset with finitely 
many components, with the components of complex codimension
$1$ being $R_1,\dots,R_s$, then it is general knot theory 
that $\pi_1(X-V;*)$ is normally generated by the classes 
of loops $l_i$, $i=l,\dots,s$, defined as in the case 
studied earlier of $X=\C^n$. In fact, even when $V$ is 
singular (without multiple components) and the $R_i$
are the complex-codimension-$1$ components of its regular 
set, the same conclusion holds---one need only observe 
that the union of the singular set $S(V)$ and the regular
components of complex codimension $2$ or more, as an 
analytic variety in its own right, has a resolution which 
is a smooth map of a smooth manifold into $X$; then any loop 
in $X-V$ may be made to bound a smooth $2$-disk in $X$ transverse 
to the resolution, and therefore disjoint from its image. Note 
however that this argument depends on the ambient space $X$ 
being a smooth manifold with its given structure as analytic space.
In this connection it is worth contemplating the example of 
$X=\{(z_1, z_2, z_3, z_4) \in \C^4\Suchthat 
z_1^2 + z_2^3 + z_3^5=0\}$. This is the product of 
$\C$ (the $z_4$ factor) with the cone on the dodecahedral 
space \cite{Milnor}, and by the celebrated Double Suspension 
Theorem\appnote{Double Suspension Theorem}, 
$X$ is homeomorpic to $\C^3$. The singular set $S(X)$ 
is a straight complex line, with real codimension $4$. 
Of course $\pi_1(X-S(X))$ has 120 elements. 
(It can be shown that each of them is, in fact, 
represented by analytic loops.)

\item 
It was asserted in the introduction that not all quasipositive 
knots were fibred.  Indeed, the first non-fibred knot in the 
Alexander-Briggs table, $5_2$, can be represented as the closure 
of the quasipositive braid $\s1^2 \s2 (\s2\s1\s2^{-1})$.

\item
For each $n$, there is an \textit{analytic} curve $V_f$ in 
$\C^2$, smooth, and $n$-sheeted over the $z$-axis, such that 
all quasipositve $n$-string closed braids occur as $K_\gamma$ 
for this $f(z,w)$ and an appropriate $\gamma$. For $n=3$, 
one may take $f(z,w)=w^3-3w + 2 \exp z$. Here, the
points of $B$ are the integral multiples of $\pi i$, 
and $B^{+}$ is a union of horizontal rays.

\item
Every oriented link has infinitely many representations 
as a closed braid (see \cite{Birman}). It would be interesting 
to have purely knot-theoretical necessary and/or sufficient 
conditions that one of the representations be quasipositive. 
Presumably not every knot or link has such a representation. 
I hope to return to this and related questions in a future paper.
\appnote{what links are qp?}
\end{inparaenum}
\end{remarks}

\noindent
{\small{\textit{Acknowledgements}---I wish to thank the referee 
of an earlier draft for pointing out that, with the conventions 
I had then established, my braids were actually quasi\textit{negative}. 
And I am very grateful to the 1980 Georgia Topology Conference, 
where in giving a talk on some of these results I first realized 
the significance of analytic disks in the configuration space.}}

\renewcommand{\MR}[1]{%
\href{http://www.ams.org/mathscinet-getitem?mr=#1}{MR#1}%
}

{}

\section*{ADDENDA}
\hypertarget{addenda}{}
Typographical errors in the original publication
have been corrected without notice; it is to be hoped
that no new ones have been introduced.  The following
notes provide updates on various points.

\def\wherewhat#1{\label{#1}}

\begin{asparaenum}
\item \wherewhat{fibration of positive braid}
\cite{Rudolph1983} gives another proof that the 
closure of a positive braid is a fibered link, by
constructing the fibration explicitly.
\item\wherewhat{Orevkov}
Some other applications of the oriented graph $B^{+}$ 
have been given by Orevkov \cite{Orevkov--Zariski}, 
\cite{Orevkov--diagrams} and Dung \cite{Dung}.
\item\wherewhat{Blank}
I am indebted to Stepan Orevkov for his observation
that Sandy Blank's unpublished 1967 thesis (see \cite{Blank})
contains a proof that (what is here called) the quasipositivity
of $\a$ is equivalent to the existence of 
an immersion $\Interior{D}''\to\C$
like that constructed in the proof of Lemma~\ref{Lemma 2}.
\item\wherewhat{Double Suspension Theorem}
The ``celebrated Double Suspension Theorem'' is expounded in \cite{DST}.
\item\wherewhat{what links are qp?}
The existence of a link which has no
representation as the closure of a quasipositive braid 
was first proved using knot polynomials, as a corollary
to an inequality of Morton \cite{Morton--inequality}
and Franks and Williams \cite{FranksWilliams}.
Boileau and Orevkov \cite{Boileau-Orevkov} 
have characterized such ``quasipositive links''
as precisely the links isotopic to boundaries of 
pieces of complex plane curve in $D^4$, but 
``purely knot-theoretical necessary and/or sufficient 
conditions'' remain elusive.
\item\wherewhat{B-W update}
\cite{Birman-Williams} was published as \cite{Birman-Williams updated}.
\item\wherewhat{Braided Surfaces update}
\cite{Rudolph1981} was published as \cite{Rudolph1981 updated}.
\end{asparaenum}

\renewcommand\refname{Additional References}

\thispagestyle{empty}


\begin{thebibliography}{10}

\bibitem{Birman}
\textsc{J\scriptsize{oan} S. B\scriptsize{irman}}: 
\textit{Braids, Links, and Mapping Class Groups}. 
Annals of Math. Studies No. 82.  Princeton University Press (1974).

\bibitem{Birman-Williams}
\textsc{J. B\scriptsize{irman}} and \textsc{R. W\scriptsize{illiams}}:
Knotted Orbits of Dynamical Systems, \textit{Topology},
to be published (1980).\appnote{B-W update}

\bibitem{Gamelin}
\textsc{T. G\scriptsize{amelin}}:
\textit{Uniform Algebras}.  Prentice--Hall, New Jersey (1969).

\bibitem{Laufer}
\textsc{H\scriptsize{enry} L\scriptsize{aufer}}:
On the number of singularities of an analytic curve.
\textit{Trans. Am. Math. Soc.} \textbf{186} (1969), 527--535.

\bibitem{Lefschetz}
\textsc{S. L\scriptsize{efschetz}}:
\textit{L'Analysis Situs et la Geometrie Algebrique.}
Gauthier--Villars, Paris (1924).

\bibitem{Milnor}
\textsc{J. M\scriptsize{ilnor}}: 
\textit{Singular Points of Complex Hypersurfaces},  
Annals of Math. Studies No. 61.  Princeton University Press (1969).

\bibitem{Rudolph1981}
\textsc{L\scriptsize{ee} R\scriptsize{udolph}}:
{Seifert Ribbons for Closed Braids}, preprint (1981).
\appnote{Braided Surfaces update}

\bibitem{Rudolph1982}
\textsc{L\scriptsize{ee} R\scriptsize{udolph}}:
{Non-Trivial Positive Braids Have Positive Signature},
\textit{Topology} \textbf{21} (1982), 325--327.

\bibitem{Stallings}
\textsc{J\scriptsize{ohn}~R.~S\scriptsize{tallings}}: 
{Constructions of fibred knots and links}, 
\textit{Proceedings of Symposia in Pure Mathematics},
Vol. XXXII, Part 2 (Providence: AMS), 1979, pp.~55--60.   

\end{thebibliography}

\begin{thebibliography}{10}
\setcounter{enumiv}{9}

\bibitem{Birman-Williams updated}
\textsc{J. B\scriptsize{irman}} and \textsc{R. W\scriptsize{illiams}}:
{Knotted periodic orbits in dynamical systems. {I}. {L}orenz's
equations}, \textit{Topology} \textbf{22} (1982), 47--82.
\MR{0682059}

\bibitem{Blank}
\textsc{V. P\scriptsize{o\'enaru}}: Extension des immersions en 
codimension $1$ (d'apr\`es Samuel Blank).  
\textit{S\'eminaire Bourbaki \textup(1967/68\textup), 
Exp. No. 342}, pp. 1--33, 
W. A. Benjamin (1969).  \MR{0255335}

\bibitem{Boileau-Orevkov}
\textsc{M\scriptsize{ichel} B\scriptsize{oileau}} and 
\textsc{S\scriptsize{tepan} Y\scriptsize{u}. O\scriptsize{revkov}}:
{Quasipositivit\'e d'une courbe analytique dans une boule pseudo-convexe},
\textit{C. R. Acad. Sci. Paris} \textbf{332} (2001), 825--830.
\MR{1836094}

\bibitem{Dung}
\textsc{N\scriptsize{guyen} V\scriptsize{iet} D\scriptsize{ung}}:
Braid monodromy of complex line arrangements.
\textit{Kodai Math. J.} \textbf{22} (1999), 46--55.
\MR{1679237}

\bibitem{FranksWilliams}
\textsc{J. F\scriptsize{ranks}} and 
\textsc{R. F. W\scriptsize{illiams}}: 
Braids and the Jones-Conway polynomial, 
\textit{Trans. Amer. Math. Soc.} \textbf{303} (1987), 97--108.
\MR{0896009}

\bibitem{DST}
\textsc{F\scriptsize{ran\,cois} L\scriptsize{atour}}:
{Double suspension d'une sph\`ere d'homologie [d'apr\`es R. Edwards]},
\textit{S\'eminaire Bourbaki, 30e ann\'ee \textup(1977/78\textup), 
Exp. No. 515}, Lecture Notes in Math. \textbf{710}, pp. 169--186, 
Springer (1979). \MR{0554220}

\bibitem{Morton--inequality}
\textsc{H. M\scriptsize{orton}}:
Seifert circles and knot polynomials, 
\textit{Math. Proc. Cambridge Philos. Soc.} \textbf{99} (1986),
107-109. \MR{0809504}

\bibitem{Orevkov--Zariski}
\textsc{S\scriptsize{tepan} Y\scriptsize{u}.\ O\scriptsize{revkov}}:
{The fundamental group of the complement of a plane algebraic curve},
\textit{Mat. Sb. (N.S.)} \textbf{137(179)} (1988), {260--270, 272}.
\MR{0971697}


\bibitem{Orevkov--diagrams}
\textsc{S\scriptsize{tepan} Y\scriptsize{u}. O\scriptsize{revkov}}:
{Rudolph diagrams and analytic realization of the {V}itushkin
covering}, 
\textit{Mat. Zametki} \textbf{60} (1996),
{206--224, 319}. \MR{1429122}


\bibitem{Rudolph1983}
\textsc{L\scriptsize{ee} R\scriptsize{udolph}}:
\href{http://arxiv.org/abs/math.GT/0106058}%
{Some knot theory of complex plane curves},
\textit{N{\oe}uds, Tresses, et Singularit\'es},
Monogr. Enseign. Math. \textbf{31} (1983), 99--122. 
\MR{0728581}

\bibitem{Rudolph1981 updated}
\textsc{L\scriptsize{ee} R\scriptsize{udolph}}:
\href{http://134.76.163.65/servlet/digbib?template=view.html&%
  id=214086&startpage=5&endpage=41&imageset-id=5267}%
{Braided surfaces and Seifert ribbons for closed braids},
  \textit{Comment. Math. Helv.} \textbf{58} (1983), 1--37. 
\MR{0699004}


\end{thebibliography}
\end{document}